\documentclass[12pt]{amsart}
\usepackage{amsmath,amssymb,amscd,amsfonts}

\newtheorem{thm}{Theorem}[section]

\newtheorem{prop}[thm]{Proposition}
\newtheorem{cor}[thm]{Corollary}

\newtheorem{assu}[thm]{Assumption}

\newcommand{\Pic}{{\text{\rm Pic}}}

\newcommand{\lra}{\longrightarrow}
\newcommand{\ot}{{\otimes}}
\newcommand{\OO}{{\mathcal O}}

\newcommand{\OX}{{\mathcal{O} _X}}

\newcommand{\inv}{^{-1}}

\newcommand{\si}{\sigma}

\newcommand{\C}{{\mathbb C}}
\newcommand{\Z}{{\mathbb Z}}

\title{Surfaces with $p_g=q=3$}
\author{Christopher D. Hacon \and Rita Pardini}
\date{}
\begin{document}

\begin{abstract}
We classify minimal complex surfaces of general type with $p_g=q=3$.
More precisely,  we show that  such a surface is either the
symmetric product of a curve of genus $3$ or a free $\Z_2-$quotient of the
product of a curve of genus $2$ and a curve of genus $3$. Our main tools are
the
generic vanishing theorems of Green and Lazarsfeld and Fourier--Mukai
transforms.
\smallskip

\noindent 2000 Mathematics Classification: 14J29.
\end{abstract}
\maketitle
\section{Introduction}
Let $X$ be a smooth minimal complex surface of general type. The smallest
possible
value of
the Euler --Poincar\'e characteristic $\chi(X)=1+p_g(X)-q(X)$ is $1$, and
for
$\chi(X)=1$ one has the bounds $1\le K^2_X\le 9$. If, in addition, the
surface
is
irregular, i.e. $q(X)>0$,  then one also has
$K^2_X\ge 2p_g(X)$  (cf. \cite{debarre}), so that $p_g(X)\le 4$. The limit
case
$p_g(X)=4$ corresponds to the product of two curves of genus $2$ (cf.
\cite{appendice}). Here we consider  the case $q(X)=3$ (cf. Theorem
\ref{main}) and
prove that
$X$ is either the symmetric product of a curve of genus  $3$ ($K^2_X=6$) or
a free
$\Z_2-$quotient of the product of a  curve of genus 3  and a curve of genus
2
($K^2_X=8$). Both surfaces have already been  described in
\cite{ccm}, where it is  also shown that the former is the  only example
with
$K^2_X=6$ and the latter is the only example with a  pencil of curves of
genus 2.

However, our approach is independent of the results of \cite{ccm}, and, as
far as we
know, it is a new one in the classification of surfaces of general type. We
are able
to  identify the surfaces by looking at  the locus $V^1(X):=\{P\in
\Pic^0(X)|h^1(-P)>0\}$, whose properties have been described very precisely
by several authors (cf. \cite{gl1}, \cite{gl2}, \cite{annulation},
\cite{simpson},
\cite{el}).
Roughly speaking, if $V^1(X)$ is $0-$dimensional we prove that  $X$ is the
symmetric
product  of a curve of genus 3 by using a cohomological characterization of
theta
divisors, due to the first author (cf. Theorem \ref{theta}). If instead
$V^1(X)$
has positive dimension, then we show the existence of a pencil of curves of
genus 2
by  using the infinitesimal description of $V^1(X)$ (cf. \cite{gl1}) and the
Castelnuovo --de Franchis Theorem.
As a byproduct of the classification, we obtain a description of the moduli
space
  of surfaces of general type with $p_g=q=3$ (cf.
Corollary 2.4) and we determine the degree
of the bicanonical map (cf. Corollary 2.3).

\paragraph{{\bf Acknowledgements}}
Gian Pietro Pirola has informed us that he has also proven
  our main result Theorem
2.2. His proof (\cite{piro}) uses different methods, based on the
techniques and results of \cite{ccm}.

 We are
indebted   to Margarida Mendes Lopes for drawing our attention to the
problem of
classifying surfaces of general type with
$p_g=q=3$,    for several useful conversations on this subject and for
pointing out
some inaccuracies in an earlier version of this paper.

This research originated during a
visit to Pisa of the first author, supported by the program ``Short term
mobility'' of
G.N.S.A.G.A. of C.N.R.
\bigskip

\paragraph{{\bf Notations and conventions}} We work over the complex
numbers; all
varieties are projective.  We use the standard notation of
algebraic geometry; we just recall here the notation for the invariants of a
surface
$X$: $K_X$ is the {\em canonical class}, $p_g(X)=h^0(X,K_X)$ the {\em
geometric
genus} and $q(X)=h^1(X,\OX)$ the {\em irregularity}. An {\em irrational
pencil of genus
$g$} on a surface $X$ is a fibration $p\colon X\to B$ with $B$ a smooth
curve of genus
$g>0$.

\section{The classification theorem}
Throughout the paper we make the following assumption:
\begin{assu} \label{ipotesi}$X$ is a smooth minimal complex projective
surface of
general
type with
$p_g(X)=q(X)=3$.
\end{assu}
We denote by $A$ the Albanese variety of $X$ and by
$a\colon X\to A$ the Albanese map.
We also assume that we have fixed a K\"ahler metric on $X$,  so that for
every
$P\in \Pic^0(X)$  and $i\ge 0$ there is an antilinear isomorphism
$H^i(X,-P)\simeq
H^0(X,\Omega^i_X(P))$ (cf. \cite{gl1}, 2.5).

Here are two examples of surfaces satisfying  Assumption
\ref{ipotesi}:

\noindent{\bf Example 1.} Let $C$ be a smooth  curve of genus $3$. Then the
symmetric product of $C$ is a surface $X$ with $K^2_X=6$ satisfying
Assumption
\ref{ipotesi}. Let  $\Pic^2(C)$ denote the subset of $\Pic(C)$ consisting of
the  line
bundles of degree 2 and let $\beta\colon X\to \Pic^2(C)$ be the map that
sends an
unordered pair
$\{p,q\}\in X$ to the linear equivalence class of $p+q$. The image $D$ of
$\beta$ is a
principal polarization by \cite{LB},  Corollary 11.2.2. By the Riemann's
Singularity
Theorem (cf. \cite{LB}, Theorem 11.2.5)  $D$ is smooth if $C$ is not
hyperelliptic,
while  it has a double point at the canonical class if $C$ is hyperelliptic.
By
Riemann--Roch on
$C$,
$\beta$ is $1-$to$-1$ if $C$ is not hyperelliptic,  while if
$C$ is
hyperelliptic  it contracts to a singular
point of type $A_1$ the $-2-$curve of $X$ corresponding to the canonical
series. In
either  case $D$ is the canonical model
of $X$. It follows that $A$ is the Jacobian of $C$ and that, up to the
choice of an
identification of $\Pic^2(C)$ with $\Pic^0(C)$, $\beta$ is the Albanese map
of $X$.
\medskip

\noindent{\bf Example 2.} Let $C_1$ be a curve of genus $2$ with an elliptic
involution $\si_1$ and $C_2$ a curve of genus $3$ with a free involution
$\si_2$.
We write $B_1:=C_1/<\si_1>$ and $B_2:=C_2/<\si_2>$. The curve $B_1$ has
genus $1$
and $B_2$ has genus $2$. We let $\Z_2$ act freely on the product $C_1\times
C_2$ via
the involution
$\si_1\times
\si_2$. The quotient surface $X:=(C_1\times C_2)/\Z_2$ is a surface with
$K^2_X=8$
satisfying Assumption  \ref{ipotesi}. The projections of $C_1\times C_2$
onto
$C_1$ and $C_2$ induce
fibrations $p_i\colon X\to B_i$, $i=1,2$. The singular fibres
of
$p_1$ are two double fibres with smooth support, occurring at the branch
points
of $C_1\to B_1$, while the fibres of $p_2$ are all smooth. The Albanese
variety
of $X$ is isogenous to the product of the Jacobians of $B_1$ and $B_2$.
\medskip

Our aim  is to prove the following:
\begin{thm}\label{main} The possibilities for a
smooth minimal surface $X$  with
$p_g(X)=q(X)=3$
 are the following:
\begin{itemize}
\item[i)] $K^2_X=6$ and $X$ is the surface of Example 1;

\item[ii)] $K^2_X=8$ and $X$ is the surface of Example 2.
\end{itemize}
\end{thm}

Before proving Theorem \ref{main} we deduce some consequences from it.
\begin{cor}
The bicanonical map     of a minimal surface of general type  $X$ with
$p_g(X)=q(X)=3$ has degree
$2$.
\end{cor}
\begin{proof}
By Theorem \ref{main} $X$ is either the surface of Example  1 or the surface
of
Example 2. In the former case the bicanonical map has degree 2 by
Proposition 3.17
of \cite{ccm}. In   the latter case, the degree is 2 by Theorem
5.6 of \cite{xiao}.
\end{proof}

\begin{cor}The moduli space  ${\mathcal M}$ of surfaces of general type with
$p_g=q=3$
has two irreducible connected components of dimensions respectively $6$ and
$5$.
\end{cor}
\begin{proof} By Theorem \ref{main}, ${\mathcal M}$ is a disjoint union
${\mathcal
M}_6\cup{\mathcal M}_8$, where ${\mathcal M}_{\alpha}:=\{[X]\in {\mathcal
M}|
K^2_X=\alpha\}$.  The sets  ${\mathcal M}_{\alpha}$ are open since $K^2$ is
a
topological invariant, thus we only have to show that ${\mathcal
M}_6$ and ${\mathcal M}_8$ are irreducible of dimension $6$ and $5$
respectively.
The points of ${\mathcal M}_6$ are in one-to-one correspondence with the
isomorphism
classes of curves of genus 3, thus the result is well known in this case.

If
$[X]\in
{\mathcal M}_8$,   then $p_2\colon X\to B_2$ is the only irrational pencil
of genus 2
of $X$. This can be seen in several ways,  for instance by observing
that a base-point free pencil is determined uniquely by the span of its
class in
$H^2(X,\C)$  and that
$H^2(X,\C)$, being of dimension 2, has only two isotropic lines,  spanned by
the
classes of the general  fibre of $p_1$ and $p_2$. It follows that the double
cover
$C_1\times C_2\to X$ is also determined uniquely, since it is the \'etale
cover of
 $X$ that kills the monodromy of the
pencil
$p_2\colon X\to B_2$. So $[X]$ is determined by the choice, up to
isomorphism,  of
a curve $C_1$ of genus $2$ with an elliptic involution $\si_1$ and of a
curve $C_2$ of
genus
$3$ with a free involution $\si_2$.
The pair $(C_1, \si_1)$ is determined by the quotient curve
$E_1:=C_1/\si_1$, by the
line bundle
$L$ of degree 1 on $E_1$ such that $\psi_*\OO_{C_1}=\OO_{E_1}\oplus L\inv$,
 where $\psi\colon C_1\to E_1$ is  the quotient map, and  by the branch
divisor
$D\equiv 2L$ of $\psi$. So, taking into account the action of the
automorphism group  of
$E_1$, we see that
$(C_1,\si_1)$ depends on 2 parametres. In addition,   it is not difficult to
write down
an irreducible  global family containing  all the  isomorphism classes of
double covers
of genus 3 of  elliptic curves. Thus  the isomorphism classes of the pairs
$(C_1, \si_1)$ form an irreducible family of dimension 2. Analogously, the
isomorphism
class  of
$(C_2,\si_2)$ is determined by the genus 2 curve $B_2:=C_2/\si_2$ and by the
choice of a
$2-$torsion line bundle $L$ of $B_2$.  The pairs $(B_2,L)$ form a
$3-$dimensional
irreducible family (cf. for instance \cite{LB}, Chapter 8, \S 3).
\end{proof}
We now turn to the proof of Theorem \ref{main}.
Surfaces satisfying Assumption \ref{ipotesi} are studied in section $3$ of
\cite{ccm}, where it is proven that if  $K^2_X=6$ then $X$ is the
surface of Example 1.   We recall some general facts from
\cite{ccm}.
\begin{prop}\label{preliminari}Let $X$ be a surface as in Assumption
\ref{ipotesi}. Then:
\begin{itemize}
\item[i)] $6\le K^2_X\le 9$;
\item[ii)]the Albanese image of $X$ is a surface;
\item[ iii)] if $X$ has an irrational pencil of genus $g>1$, then $X$ is the
surface of
Example 2.
\end{itemize}
\end{prop}
\begin{proof}Statement i) follows from Miyaoka--Yau inequality  $K^2_X\le
9\chi$ and from the inequality  $K^2_X\ge 2p_g$ of \cite{debarre}.
Statements ii)
and iii) correspond  to Proposition 3.1, i) and Theorem 3.23 of \cite{ccm},
respectively.
\end{proof}
The proof of Theorem \ref{main} is based on the study of the set
$V^1(X):=\{P\in
\Pic^0(X)|\ h^1(-P)>0\}$. The  sets
$V^i(Y):=\{P\in
\Pic^0(Y)|\ h^i(-P)>0\}$, for $Y$ a  variety of any dimension,  have  been
studied  by
Green
and Lazarsfeld (\cite{gl1}, \cite{gl2}), by Simpson (\cite{simpson}), and in
the
case of surfaces  by Beauville  (\cite{annulation}). We  recall here only
those
properties of
$V^1(X)$ that we are going to use:
\begin{thm}\label{V1varie}Let $X$ be a complex surface of maximal Albanese
dimension. Then:
\begin{itemize}
\item[i)] $V^1(X)$ is a proper subvariety of $\Pic^0(X)$;
\item[ii)]the irreducible
components of
$V^1(X)$ are translates of abelian subvarieties of $\Pic^0(X)$ by torsion
points;
\item[iii)] $V^1(X)$ has a component of dimension $>1$ iff $X$ has an
irrational pencil
of genus $>1$:
\item[iv)] Let $P\in V^1(X)$ be a point, let  $v\in H^1(X,\OO_X)$ be a
nonzero
 vector and let
$\si\in H^0(X,\Omega^1_X)$ be the conjugate of $v$. If $v$ is in the tangent
cone to
$V^1(X)$ at $P$, then the map $H^0(X, \Omega^1_X(P))\stackrel{\wedge
\si}{\lra}
H^0(X,
\omega_X(P))$ is not injective;
\item[v)]Let  $P\in V^1(X)$. Then $\{P\}$ is a component of $V^1(X)$ iff
the map $H^0(X, \Omega^1_X(P))\stackrel{\wedge \si}{\lra} H^0(X,
\omega_X(P))$ is injective for all nonzero $\si\in H^0(X,\Omega^1_X)$.
\end{itemize}
\end{thm}
\begin{proof}
Statement i) is Theorem 1 of \cite{gl1}.
The fact that the components of $V^1(X)$ are translates of  abelian
subvarieties
follows from Theorem 0.1 of \cite{gl2} and the fact  that the translation is
by a
torsion point follows from Theorem 4.2 of \cite{simpson}. Statement iii)
follows from
\cite{annulation}, Corollary 2.3.  Statement iv) follows by combining
Theorem
1.6, Lemma 2.3 and Lemma 2.6 of \cite{gl1}. Statement v) is a consequence of
Theorem 1.2,
(1.2.3) of \cite{el}.
\end{proof}
If $X$ is the surface of Example 2, then the set  $V^1(X)$ has
a component of dimension 2 by Theorem \ref{V1varie}, iii).  On the other
hand, using
adjunction on
$A:=\Pic^2(C)$, it is easy to show that for the surface of Example 1 one has
$V^1(X)=\{\OO_X\}$.
The results that follow show that the structure of $V^1(X)$ characterizes
$X$.

\begin{prop}\label{dimzero}
If $X$ does not have an irrational pencil of genus $>1$, then the
dimension of
$V^1(X)$ is
$0$.
\end{prop}
\begin{proof}
 We remark first of all that, in view of the assumption, we have    $\dim
V^1(X)\le
1$ by Theorem \ref{V1varie}, iii).

Assume now that there is a $1-$dimensional component $T$ of $V^1(X)$. By
Theorem \ref{V1varie}, ii),
$T= T_0+Q$, where $T_0$ is an abelian subvariety of $\Pic^0(X)$ and $Q$ is a
torsion point. Notice that $Q\notin T_0$, since by Proposition 4.1 of
\cite{gl1}
$\OO_X$ is an isolated point of $V^1(X)$. Fix a torsion point $P\in T$,
denote by
$n$ its order and assume that
$n$ is minimal. This is the same as saying that if $kP\in T_0$ then
$kP=0$. Let
$v\in H^1(\OO_X)$ be a nonzero vector tangent  to $T$ and  let $\si\in
H^0(\Omega^1_X)$ be the conjugate of $v$. The vector $v$ lies
in the tangent cone to $V^1(X)$
at $P$ and therefore the map $H^0(\Omega_X^1(P))\stackrel{\wedge \si}{\to}
H^0(\omega_X(P))$ is not injective by Theorem \ref{V1varie}, iv). Denote by
$\tau$ a nonzero element in the kernel of this map. Let
$\pi\colon Y\to X$ be the connected
\'etale cyclic cover associated to the finite subgroup  $<\! P\!>\subset
\Pic^0(X)$. The group $<\!P\!>$ can be naturally  identified with the dual
group of the Galois group $G\cong\Z_n$  of $\pi$,  and
$\pi^*H^0(\Omega^1_X(P))$
is the eigenspace of $H^0(\Omega^1_Y)$ on which $G$ acts via the character
corresponding to
$P$. Since $\pi^*\si\wedge \pi^*\tau=0$, by
the Castelnuovo--de Franchis Theorem there
exists a fibration $q\colon Y\to B$ onto a curve
of genus at least $2$,   such
that $\pi^*\si$ and $\pi^*\tau$ are
pull-backs from $B$. The fibres of $q$ are integral curves of $\pi^*\si$.
Since $\pi^*\si$ is $G-$invariant, it follows that $G$ permutes the fibres
of $q$ and thus induces a pencil  $p\colon X\to E:=B/G$. More precisely,  we
have a commutative diagram:
\[
\begin{CD}
Y@>\pi>>X\\
@VqVV@VpVV\\
B@>\bar{\pi}>>E
\end{CD}
\]
 Let $\bar{\si}$ be the $1-$form on $B$
such that $\pi^*\si=q^*\bar{\si}$. By commutativity of the diagram,
$\bar{\si}$ is
$G-$invariant and thus it is a pull-back from $E$. This shows that $E$ is
not
rational. In view of the assumptions, we conclude that $E$ has genus $1$.
Since
$\si$ is a pull-back from
$E$, comparing the tangent spaces at the origin one sees that
$T_0=p^*\Pic^0(E)$.

 Now let $\bar{\tau}$ be the $1-$form on $B$ such that
$\pi^*\tau=q^*\bar{\tau}$. By commutativity of the diagram, $\bar{\tau}$ is
an
eigenvector for $G$ with character corresponding to $P$. Since this
character
generates the dual group $G^*$, the action of $G$ on $B$ is effective, i.e.
$\bar{\pi}\colon B\to E$ is a $G-$cover.

We claim that $\bar{\pi}$ is totally ramified, namely it does not factor
through
an \'etale cover of $E$. Assume otherwise, and denote by $H\subsetneq G$ the
subgroup generated by the elements that do not act freely on $B$. Set
$Y':=Y/H$,
$B':=B/H$. Then we have a commutative diagram:
\[
\begin{CD}
Y'@>\pi'>>X\\
@Vq'VV@VpVV\\
B'@>\bar{\pi}'>>E
\end{CD}
\]
where the maps are the obvious ones. In particular, $\bar{\pi}'$ is \'etale
by
construction. If we
denote by $d<n$ the order of $H$, then the cover $\pi'$
corresponds to the subgroup of $\Pic^0(X)$ generated by $dP$. On the other
hand, $\pi'$ is obtained from $\bar{\pi}'$ by taking base change with $p$,
so the
subgroup of $\Pic^0(X)$ corresponding to $\pi'$ is actually a subgroup of
$p^*\Pic^0(E)$, i.e. $dP\in p^*\Pic^0(E)=T_0$. This contradicts the
minimality
assumption
on $n$, and we conclude that $\bar{\pi}$ is totally ramified.
The $G-$action on $B$ gives a decomposition
$\bar{\pi}_*\OO_B=\OO_E\oplus_{\chi\in G^*-\{1\}}L_{\chi}\inv$, where the
$L_{\chi}$ are line bundles and $L_{\chi}\inv$ is the eigenspace of
$\bar{\pi}_*\OO_B$ on which $G$ acts via the character $\chi\in G^*$. Saying
that
$\bar{\pi}$ is totally ramified is the same as saying that $\deg L_{\chi}>0$
for
every $\chi\in G^*-\{1\}$. The standard formulas for abelian covers (cf.
\cite{ritaabel}, Proposition 4.1), give $g(B)=g(E)+\sum_{\chi\in
G^*-\{1\}}h^0(\omega_E\ot L_{\chi})=1+\sum_{\chi\in G^*-\{1\}}\deg
L_{\chi}\ge n$.
Denote by $f$ the genus of the general fibre of $q$. The  Corollary on page
344 of
\cite{appendice} gives:
$\chi(Y)=n\chi(X)=n\ge (f-1)(g(B)-1)\ge (f-1)(n-1)$. Thus we either have
$f=2$ or
$n=2$, $f=3$. Recall that $f$ is also the genus of the general fibre of $p$.
Then $f=2$ is impossible: indeed by the Lemma on page 345 of
\cite{appendice} the equality
$3=q(X)=g(E)+f$ would imply that
$X$ is the product of $E$ with a curve of genus $2$, contradicting the fact
that $X$
is of general type. So we are left with the case $n=2$, $f=3$. Notice that
in this
case we also have
$g(B)=2$, so that $2=\chi(Y)=(g(B)-1)(f-1)$. So, again by the  Corollary of
page 344
of
\cite{appendice}, the fibration $q$ is isotrivial with smooth fibres, and
the only
singular fibres of $p$ are two double fibres with smooth support, occurring
at the two
branch points of $\bar{\pi}\colon B\to E$.    Hence by \cite{serrano},
Theorem 2.1 and
Remark 2.3, there exist a curve
$C_1$ of genus
$3$, a curve
$C_2$  of genus
$g_2$, and a finite group
$G$ acting faithfully  both on
$C_1$ and on $C_2$ such that:

i) $C_2/G$ is isomorphic to $E$;

ii) the diagonal action on $C_1\times C_2$ defined by
$(x,y)\stackrel{g}{\mapsto}(gx,gy)$ is free;

iii) $X$ is isomorphic to the quotient $(C_1\times C_2)/G$, and the
fibration
$p$ corresponds to the map $(C_1\times C_2)/G\to C_2/G=E$ induced by the
second
projection $C_1\times C_2\to C_2$.

In addition, the following hold:

a) $C_1/G$ has genus $2$ (since $q(X)=g(C_2/G)+g(C_1/G)=1+g(C_1/G)$ by
Proposition 2.2
of \cite{serrano});

b) $2(g_2-1)=|G|$ (since $2(g_2-1)=\chi(C_1\times C_2)=|G|\chi(X)=|G|$).

Applying the Hurwitz formula to the quotient map $C_1\to C_1/G$ we see that
$G$ has
order $2$. Condition b) now implies  $g_2=2$, hence $X$ is
the surface of Example 2 and $p$ is the
pencil $p_1$. Thus $X$   has also an irrational pencil of genus $2$,
contrary to the
assumptions.
\end{proof}
\begin{prop}\label{V1OX}If $X$ does not have an irrational pencil of genus
$>1$,
then $V^1(X)=\{\OO _X\}$.
\end{prop}
\begin{proof} Let $\OO _X\ne P\in V^1(X)$. By Proposition \ref{dimzero},
$\{P\}$ is a
component of $V^1(X)$. It follows from Theorem \ref{V1varie}, v) that for
all $0\ne
\si
\in H^0(X,\Omega^1_X)$ the map $H^0(X,\Omega^1_X(P))\stackrel{\wedge
\si}{\lra}
H^0(X,\omega_X(P))$  is injective. Using Hodge theory with twisted
coefficients  it
follows that the map  induced by cup product
$$H^1(-P) \otimes H^1(X,\OO _X) \longrightarrow H^2(-P)$$
is non zero on (non zero) simple tensors. Therefore, by a result of H. Hopf
(see \cite{acgh} pg. 108)
$h^2(-P)\geq h^1(-P)+h^1(\OO _X)-1$. Hence
$$1=\chi (\OO _X)=\chi (-P)=h^2(-P)-h^1(-P)\geq 2.$$
This is the required contradiction.
\end{proof}
The last ingredient of the proof of Theorem \ref{main} is the following
result from
\cite{fourier}:
\begin{thm}\label{theta}
Let $Z$ be a smooth complex projective variety of dimension $n$, let $A$ be
an abelian
variety and let
$f\colon Z\to A$ be a generically finite morphism such that $f(Z)$ is a
divisor.
Assume that:

i) $h^0(Z, \Omega^i_Z)=$ ${n+1}\choose{i}$ for $0\le i\le n$;

 ii) $h^i(Z,
\omega_Z(f^*P))=0$ for all $\OO_A\ne P\in \Pic^0(A)$ and all $i>0$.

\noindent Then $A$ is
principally polarized and $f(Z)$ is a theta divisor.
\end{thm}
\begin{proof}
Cf. \cite{fourier}, Corollary 3.4.
\end{proof}
\begin{proof}[Proof of Theorem \ref{main}]
By Proposition \ref{preliminari}, iii), $X$ has an irrational pencil of
genus $>1$
iff
$X$ is the surface of Example 1.

If $X$ has no irrational pencil of genus $>1$, then
$V^1(X)=\{\OO_X\}$
by Proposition \ref{V1OX}. Thus we can apply Theorem \ref{theta} to
the Albanese map $a\colon X\to A$. It follows that $A$ is pricipally
polarized
and
$Y:=a(X)$ is a  theta divisor. It is well known (cf.
\cite{LB}, Ch. 11,  Corollary 8.2,b)) that an abelian
 threefold with an irreducible principal polarization is the Jacobian of a
curve $C$ of genus $3$. As already explained in Example 1, the theta divisor
is
the canonical model of the symmetric product of $C$. Thus $K_Y$ is an ample
Cartier divisor with $K^2_Y=6$. To finish the proof, it suffices to show
that the
map
$a\colon  X
\to Y $ is birational.  If we denote by
$R$ the divisorial part of the ramification locus of $a$, the adjunction
formula gives $K_X=a^*K_Y+R$. If we denote by $d$ the degree of $a$,
then we have:
$$K^2_X=a^*K_YK_X+RK_X\ge a^*K_YK_X\ge (a^*K_Y)^2=6d$$
where the first inequality follows from the fact that $K_X$ is nef and the
second one from the fact that $a^*K_Y$ is nef. Since $K^2_X\le 9$ by
Proposition \ref{preliminari}, i), it follows that $d=1$ and $X$ and $Y$ are
isomorphic.
\end{proof}

\bigskip
\bigskip

\begin{minipage}{13cm}
\parbox[t]{5.5cm}{Christopher D. Hacon\\
Department of Mathematics\\
Sproul Hall 2208\\
University of California\\
Riverside, CA 92521-013 USA\\
hacon@math.ucr.edu}
\hfill
\parbox[t]{5.5cm}{Rita Pardini\\
Dipartimento di Matematica\\
Universit\`a di Pisa\\
Via Buonarroti, 2\\
56127 Pisa, ITALY\\
pardini@dm.unipi.it}
\end{minipage}

\begin{thebibliography}{ACGH}
\bibitem[ACGH]{acgh} E. Arbarello, M. Cornalba, P. Griffiths, J. Harris,
{\em
Geometry of algebraic curves, volume I}, G. m.
W. {\bf 267}, Springer Verlag 1985.
\bibitem[Be1]{annulation} A. Beauville, {\em Annulation du $H^1$ pour les
fibr\'es en droites plats},
Complex Algebraic Varieties, Proc. Conf., Bayreuth/Ger. 1990, Springer
L. N. M. {\bf 1507}
(1992), 1--15.
\bibitem[Be2]{appendice} A. Beauville, {\em L'inegalit\'e $p_g\ge 2q-4$ pour
les
surfaces de type g\'en\'eral}, Appendix to \cite{debarre}, Bull. Soc. Math.
de France, vol. {\bf 110}  (1982),
343--346.
\bibitem[CCM]{ccm} F. Catanese, C. Ciliberto, M. Mendes Lopes, {\em On the
classification
of irregular surfaces of general type with non birational bicanonical map},
Tran.  A.M.S., vol. {\bf 350} (1998), 275--308.

\bibitem[De] {debarre} O. Debarre, {\it  In\'egalit\'es num\'eriques pour
les surfaces de
type g\'en\'eral}, Bull. Soc. Math. de France, vol. {\bf 110}  (1982),
319--342.

\bibitem[EL] {el}
 L. Ein, R. Lazarsfeld,
{\em Singularities of theta divisors, and birational geometry of
irregular varieties},
 Jour. A.M.S.
  {\bf 10},  1 (1997),  243--258.

\bibitem[GL1]{gl1}
M. Green, R. Lazarsfeld,
{\em Deformation  theory, generic vanishing theorems,
and some conjectures of Enriques, Catanese and Beauville},
 Invent. Math.
 {\bf 90} (1987 ),  389--407.

\bibitem[GL2]{gl2}
M. Green, R. Lazarsfeld,
{\em Higher obstruction to deforming cohomology groups of line bundles},
 Jour. A. M. S.
  {\bf 4} (1991),  87--103.

\bibitem[Hac]{fourier}
C. Hacon,
{\em Fourier transforms, generic vanishing theorems and polarizations
of abelian varieties}, Math. Zeit. (to appear,
 Preprint math.AG/9902078 1999).

\bibitem[LB]{LB} H. Lange, C. Birkenhake
{\em Complex abelian varieties}, G. m. W. {\bf 302},
Springer-Verlag, 1992.

\bibitem[Pa]{ritaabel}
R. Pardini,
{\em Abelian covers of algebraic
varieties},
 J. reine angew. Math. {\bf 417} (1991), 191--213.
\bibitem[Pi]{piro} G. P. Pirola, {\it Algebraic surfaces with $p_g=q=3$ and
no
irrational pencils}, preprint.
\bibitem[Se] {serrano} F. Serrano, {\it Isotrivial fibred surfaces}, Annali
di
Matematica pura ed applicata (IV), {\bf CLXXI} (1996),  63--81.

\bibitem[Si]{simpson}
C. Simpson,
{\em Subspaces of moduli spaces of rank one local systems},
 Ann. Sci. \'{E}cole Norm. Sup. 4), {\bf 26}, 3,  (1993),
361--401.

\bibitem[Xi]{xiao} G. Xiao, {\em Surfaces fibr\'ees en courbes de genre
deux}, Springer
L. N. M.  {\bf 1137} (1985).
\end{thebibliography}
\end{document}